\documentclass[12pt]{article}

\usepackage{amssymb}
\usepackage{latexsym}
\usepackage{a4wide}
\usepackage{amscd}
\usepackage{graphics}
\usepackage{amsmath}
\usepackage{bbold}

\newcommand{\N}{\mathbb{N}}                     
\newcommand{\Z}{\mathbb{Z}}                     
\newcommand{\R}{\mathbb{R}}                     
\newcommand{\set}[2]{\left\{{#1}\mid{#2}\right\}}       
\newcommand{\proof}{{\sl Proof.}\hspace{5pt}}   
\newcommand{\qed}{\hfill $\Box$ \bigskip}       
\newcommand{\lip}{\mathrm{lip\,}}       
\newcommand{\ran}{\mathrm{ran\,}}      

\newtheorem{thm}{Theorem}[section]      
\newtheorem{introthm}{Theorem}

\newtheorem{cor}[thm]{Corollary}        
\newtheorem{lem}[thm]{Lemma}            
\newtheorem{prop}[thm]{Proposition}     
\newtheorem{defn}[thm]{Definition}
\newtheorem{rem}[thm]{Remark}

\title{On the global stable manifold}
\author{Alberto Abbondandolo, Pietro Majer}

\begin{document}

\renewcommand{\theenumi}{\roman{enumi}}
\renewcommand{\labelenumi}{(\theenumi)}

\maketitle

\begin{abstract}
We give an alternative proof of the stable manifold theorem as an
application of the (right and left) inverse mapping theorem on a
space of sequences. We investigate the diffeomorphism class of the
global stable manifold, a problem which in the general Banach setting
gives rise to subtle questions about the possibility of extending germs of
diffeomorphisms. 
\end{abstract}

\bigskip

\centerline{AMS Subject Classification: 37D10, 46B20, 58B10.}

\bigskip

Let $M$ be an open subset of a manifold $M^{\prime}$
of class $C^k$, $1\leq k\leq \omega$ (meaning that $k\in \Z^+$, or
$k=\infty$, or $k=\omega$, where as usual $C^{\omega}$ denotes the 
analytic category) modeled on the Banach space $E$. 
Let $f:M\rightarrow f(M)\subset M^{\prime}$ 
be a $C^k$ diffeomorphism\footnote{Usually $M^{\prime}=M=f(M)$, so that $f$
  is a diffeomorphism of $M$. We allow this slightly more general setting
  to include the case of a local diffeomorphism.}, 
and let $x\in M$ be a hyperbolic fixed point of $f$. This means that 
the spectrum of $Df(x)$ does not meet the unit circle, thus it is
divided into two disjoint closed subsets $\sigma(Df(x))\cap \{|z|<1\}$ and
$\sigma(Df(x))\cap \{|z|>1\}$, and the spectral decomposition theorem
gives a corresponding $Df(x)$-invariant decomposition of 
the tangent space of $M$ at $x$, $T_x M = E^s \oplus E^u$.
The stable manifold of $x$ is the set
\[
W^s(x) = \set{p\in \bigcap_{n\in \N} f^{-n}(M)}{
\lim_{n\rightarrow \infty} f^n(p)= x}.
\]
The stable manifold theorem states that $W^s(x)$ is an immersed $C^k$
submanifold of $M$. A first way to prove such a result is
to define the local stable manifold near $x$, to use the graph
transform method to show that in local coordinates such a set is the
graph of a Lipschitz map, then to prove further
regularity, and finally to use the map $f$ to describe the whole
stable manifold.
See, for example, \cite{shu87}. Another approach is the
proof by Irwin \cite{irw70}, 
where one replaces the graph transform method by 
an argument involving the implicit mapping theorem applied to the
space of sequences (see also \cite{wel76}).  
See \cite{hps77,irw80,wig94,cfl03,cfl03b} for many generalizations and 
for more bibliography. 

The first aim of this note is to give a different proof of the stable
manifold theorem. Denote by $c_x(M)$ the set of
all $M$-valued sequences converging to $x$, and recall that $c_x(M)$
has a natural structure of $C^k$ manifold modeled on the Banach space
$c_0(E)$, the space of infinitesimal $E$-valued sequences with the
supremum norm (see section \ref{uno}). The stable manifold theorem can
then be stated in the following form: 
  
\begin{introthm}
\label{stab}
Let $M$ be an open subset of the $C^k$ Banach manifold $M^{\prime}$, let
$f:M \rightarrow f(M)\subset M^{\prime}$ be a $C^k$ diffeomorphism,
$1\leq k\leq \omega$, and let $x\in M$ be a hyperbolic fixed
point of $f$, inducing the $Df(x)$-invariant splitting $T_x M= E^s
\oplus E^u$. Then the set
\[
\mathcal{W} = \set{u\in c_x(M)}{u(n+1)=f(u(n)) \; \forall n\in \N}
\]
is a closed $C^k$ submanifold of $c_x(M)$, and the evaluation map at zero,
\[
\mathrm{ev}_0: c_x(M) \rightarrow M, \quad u\mapsto u(0),
\]
subordinates a $C^k$ injective immersion of $\mathcal{W}$ 
onto $W^s(x)$, such that $D
\mathrm{ev}_0 (x) T_x \mathcal{W} = E^s$ (here $x$ denotes the
constant sequence $x(n)=x$). 
\end{introthm}

Therefore $W^s(x)$ is a $C^k$ immersed submanifold of $M$. 
Simple examples show
that $W^s(x)$ may not be an embedded submanifold, and may not be
even locally closed. Lifting the dynamical system to the space of
sequences produces instead the closed embedded submanifold $\mathcal{W}$. 
Notice also that $\mathrm{ev}_0|_{\mathcal{W}}$ is a semi-conjugacy between
the shift operator on $\mathcal{W}$ and the restriction of $f$ to
$W^s(x)$. By considering $f^{-1}$ instead of $f$ one finds an analogous
statement for the unstable manifold $W^u(x)$.

As in Irwin's approach, the proof of Theorem \ref{stab} uses the 
implicit mapping theorem on the space of sequences, but in a more 
direct way, due to the fact that we deal with zero sets of
mappings instead of graphs. The whole analysis is reduced to quite a
simple linear problem, and the regularity in the $C^k$ or even 
in the analytic case follows directly. 

Let us sketch the proof of Theorem \ref{stab} in the case $M^{\prime}
= M = \R^n$.  Assuming that the fixed point $x$ is the origin,
$\mathcal{W}$ is the zero set of the $C^k$ map
\[
F: c_0(\R^n) \rightarrow c_0(\R^n), \quad u \mapsto S u - f_*(u),
\]
where $S$ is the left shift $(Su)(n)=u(n+1)$, and $f_*$ denotes
the composition by $f$. It is readily seen that the set of $u\in
\mathcal{W}$ such that $DF(u)$ has a right inverse is open and
$f_*^{-1}$-invariant. Such a set contains the origin, because the
linear mapping
\[
DF(0) : c_0(\R^n) \rightarrow c_0(\R^n), \quad v \mapsto S v - Df(0)_*
v,
\]
has a right inverse in the form of a convolution operator (see Lemma
\ref{1} below). The only $f_*^{-1}$-invariant neighborhood of $0$ in
$\mathcal{W}$ is the whole $\mathcal{W}$, so the implicit mapping
theorem implies that $\mathcal{W}$ is a $C^k$ submanifold of
$c_0(\R^n)$. The restriction to $\mathcal{W}$ of the 
map $\mathrm{ev}_0$ is clearly injective, and an argument similar to 
the one used above shows that it is an immersion.

The usual statement of the local stable manifold theorem is then
deduced from Theorem \ref{stab} as a simple corollary (see section
\ref{tre}). The same idea works for continuous-time dynamical
systems, i.e.\ flows obtained by integrating some vector field on
$M$ (see Remark \ref{rof}).

Finally, we investigate the diffeomorphism class of the
global stable manifold $W^s(x)$, when $f$ is a global diffeomorphism 
of the Banach manifold $M$. In this case, it is natural to expect 
$\mathcal{W}$ to be $C^k$ diffeomorphic to the Banach space $E^s$, so 
that the stable manifold is a $C^k$-immersed copy of $E^s$. 
We do not know if this is true for manifolds $M$ modeled on an
arbitrary Banach space. The
difficulty in proving such a result is due to the fact that on Banach
spaces the problem of extending the germ of a map by keeping the same
regularity is quite delicate, because there need not
exist a smooth norm, or smooth partitions of unity. 
We characterize those diffeomorphisms for which $\mathcal{W}$ is
diffeomorphic to $E^s$ in terms of an extension property for the germ
of $f$ at $x$ (see Corollary \ref{dopo}), and we deduce the
following result:

\begin{introthm}
\label{diffcl}
Let $f$ be a $C^k$ diffeomorphism of the $C^k$ Banach manifold $M$,
$1\leq k \leq \omega$, and let $x\in M$ be a hyperbolic fixed point 
of $f$, with associated splitting $T_x M = E^s \oplus E^u$. Then the
manifold $\mathcal{W}$ is homeomorphic to the Banach space $E^s$, by a
bi-locally Lipschitz homeomorphism\footnote{A map $h$ between Banach
  manifolds (on which there is no preferred metric) is said to be locally
  Lipschitz if $\varphi \circ h \circ \psi^{-1}$ is locally Lipschitz,
  for every pair of local charts $\varphi$ and $\psi$. A homeomorphism
  $h$ is said to be bi-locally Lipschitz if both $h$ and $h^{-1}$ are locally
  Lipschitz.}. Assume moreover that the  
Banach space $E^s$ has the following property: there exists a {\em
  bounded} $C^k$ map $\varphi: E^s \rightarrow E^s$ such that
\begin{enumerate}
\item $\varphi$ coincides with the identity in
  a neighborhood of $0$, in the case $1\leq k <\infty$;
\item $\varphi(0)=0$ and $D\varphi(0)=I$, in the case $k=\infty$ or 
  $k=\omega$.
\end{enumerate}
Then $\mathcal{W}$ is $C^k$ diffeomorphic to the Banach space
$E^s$. 
\end{introthm}

Therefore $W^s(x)$ is always the image of $E^s$ by a locally Lipschitz
and locally closed injective map, which can be chosen to be a $C^k$ immersion
whenever the Banach space $E^s$ has one of the properties described above.

Notice that these properties are hereditary, in the sense that if
a Banach space $E$ admits a map $\varphi$ with one of these
properties, then every complemented linear
subspace of $E$ does 
(if the subspace $E^{\prime}$ is the image of the bounded linear projection
$P$, consider $P\varphi|_{E^{\prime}}$).
 
The property of admitting a $C^k$ bounded map coinciding with the
identity in a neighborhood of $0$ was introduced by Atkin
\cite{atk01}. He observed that this property not only holds trivially when
the Banach space $E$ admits a $C^k$ norm (so for instance when $E$ is
a Hilbert space), but it holds also for some non-regular Banach spaces such
as $\ell_{\infty}$ or $C^0(K,\R)$. This property implies that $C^k$ germs
at $0\in E$ have $C^k$ extensions to the whole Banach space $E$, a fact which
is useful in order to make global constructions, also in lack of
$C^k$ partitions of unity. 

Clearly, an Atkin map cannot exist in the analytic category, but in
this case (and actually also in the smooth category) it is enough to
assume the weaker condition (ii). Again, some non-regular Banach
spaces, such as $\ell_{\infty}$ and $C^0(K,\R)$, admit an analytic map 
satisfying (ii).

We do not know whether there exist Banach spaces which do not admit a
$C^k$ map satisfying (i) or (ii).    

\section{Notations, Definitions and Basic Facts}
\label{uno}

\paragraph{Linear Operators and splittings.} 
\label{upu}
Let $(E, |\cdot|_E)$ and $(F,|\cdot|_F)$ be Banach spaces. 
We denote by $\mathcal{L}(E,F)$ the Banach space of all linear bounded 
operators from $E$ to $F$, endowed with the operator norm 
$\|T\|:=\sup_{|x|_E\leq 1}|Tx|_F$. 
If $E=F$ we simply write $\mathcal{L}(E)$ for $\mathcal{L}(E,E)$. 
A linear subspace $X$ (necessarily closed) of $E$ {\em splits} if and only 
if there exists a subspace $Y$ such that $E=X\oplus Y$. 
If $L\in \mathcal{L}(E,F)$ and $R\in \mathcal{L}(F,E)$ are such that 
$LR=1_F$, $L$ is called a {\em left inverse} of $R$ or a 
{\em linear retraction} and $R$ is called a {\em right inverse} of $L$ 
or a {\em linear section}. 
Then $L$ is surjective, $R$ is injective and $E$ decomposes as 
$E=\ker L\oplus\ran R$, with projections 
$P_{\ran R} = RL$ and $P_{\ker L}=1_E - RL$.
Conversely, if  $R\in \mathcal{L}(F,E)$ is injective and
$E=X\oplus\ran R$ for some subspace $X$ of $E$, then $L:=R^{-1}P_{\ran
R}\in \mathcal{L}(E,F)$ is a right inverse of $R$, with $\ker L = X$. 
Similarly, if $L\in \mathcal{L}(E,F)$ is surjective and $E=\ker L\oplus Y$ for
some subspace $Y$ of $E$, then $R:=(L_{|Y})^{-1}$ is a
right inverse of $L$ with $\ran R=Y$. The set of linear sections $L\in
\mathcal{L}(E,F)$ and the set of linear retractions $R\in
\mathcal{L}(E,F)$ are open in $\mathcal{L}(E,F)$.
 
\paragraph{Immersions and submersions.} 
Let $M$, $N$ be differentiable manifolds of class $C^k$, 
$1\leq k\leq\omega$, modeled on the Banach space $E$, 
respectively $F$. A map $f:M\rightarrow N$ is a {\em local immersion} 
(resp.\ a {\em local submersion}) {\em at $p$}, if $f$ is a linear section 
(resp.\ a linear retraction) in local charts at $p$, meaning
that there exist a local chart at $p$, $\varphi:U\rightarrow\varphi(U)
\subset E$,  a local chart at $q:=f(p)$,
$\psi:V\rightarrow\psi(V)\subset F$, 
and a linear operator $A\in \mathcal{L}(E,F)$ which is a linear
section (resp.\ a linear retraction) and such that 
$\psi f\varphi^{-1}=A_{|\varphi(U)}$. Then $Y:=f(U)$ is 
submanifold of $N$ and its tangent space at $q$ is $T_q Y=\ran Df(p)$
(resp.\  $X:=f^{-1}(q)\cap U$ is a 
submanifold of $M$ and its tangent space at $p$ is $T_p X =\ker Df(p)$ ).
The map $f$ is said to be simply an {\em immersion} (resp.\ a {\em
submersion}), if it is a local immersion (resp.\ a local submersion) at any
$p\in M$. 
In the first case, if $f$ is also injective, $f(M)$ is
said to be an {\em immersed submanifold} of $N$. 
If $f$ is a local submersion at every $p\in f^{-1}(q)$, then
$f^{-1}(q)$ is an embedded submanifold of $M$.

The implicit mapping theorem implies
the usual criterion for local immersions and submersions,  
stating that $f$ is a local immersion (resp.\ a
local submersion) at $p$ if and only if $Df(p)\in \mathcal{L}(T_p M , T_q N)$
is a linear section (resp.\ a linear retraction). A standard reference is 
\cite{lan99}, section II, \S 2.

The criterion for local submersions has the following immediate consequence:

\begin{prop}
\label{p1}
Let $f,g:M\rightarrow N$ be $C^k$ maps between $C^k$ Banach manifolds,
$1\leq k\leq \omega$, and set
\[
W=\set{p\in M}{f(p)=g(p)}.
\]
If for every $p\in W$, the operator $Df(p)-Dg(p)\in \mathcal{L}(T_p
M,T_{f(p)} N)$ is a linear retraction, then $W$ is a $C^k$ submanifold
of $M$, with $T_p W = \ker (Df(p)-Dg(p))$.
\end{prop}

Indeed, the matter being local, we may assume that $N$ is an open subset
of the Banach space $F$, so that $W$ is the zero set of the map
$f-g$, which is by hypothesis a local submersion at every $p\in W$.

\paragraph{Discrete convolutions on $\mathbf{\ell_p}$ classes.} 
If $(E,|\cdot|)$ is a Banach space, the $\ell_p$-norm of 
$u:\Z\rightarrow E$ is $\|u\|_p:=\left(\sum_{n\in\Z}|u(n)|^p
\right)^\frac{1}{p}$, for $1\leq p<\infty$, or 
$\|u\|_\infty:=\sup_{n\in\Z}|u(n)|$.  
Then $\ell_p(\Z,E)$ denotes the Banach space of all $u:\Z\rightarrow E$
such that $\|u\|_p<\infty$. The set 
\[
c_0(\Z,E):=\set{u:\Z\rightarrow E}{\lim_{|n|\rightarrow \infty}
  u(n)=0}
\]
is a closed subspace of $\ell_\infty(\Z,E)$ and for all 
$1\leq p\leq q<\infty$, $\ell_p(\Z,E)\subset \ell_q(\Z,E)
\subset c_0(\Z,E)\subset \ell_\infty(\Z,E)$.
The analogous class $\set{u:\N\rightarrow
  E}{\lim_{n\rightarrow \infty} u(n)=0}$
is denoted simply by $c_0(E)$; it can be viewed as a closed splitting subspace 
of $c_0(\Z,E)$. Indeed, the identity mapping on $c_0(E)$ factors as 
$c_0(E)\stackrel{j}{\rightarrow}c_0(\Z,E) \stackrel{\rho}{\rightarrow}c_0(E)$,
the inclusion $j$ being given by zero-extension, the map $\rho$ being
given by restriction to $\N\subset\Z$.
 
If $g\in \ell_1(\Z,\mathcal{L}(E))$ and 
$u\in \ell_\infty(\Z,E)$, their {\em convolution product} 
$g*u$ is defined by 
\[
(g*u)(n):=\sum_{h\in\Z}g(n-h)u(h).
\]
Young's inequality $\|g*u\|_p\leq \|g\|_1 \|u\|_p$ implies that for any 
$p\in [1,+\infty]$ the convolution product is continuous
as a bilinear map $\ell_1(\Z, \mathcal{L}(E))\times 
\ell_p(\Z, E)\rightarrow \ell_p(\Z,E)$.
Furthermore, $g*u\in c_0(\Z,E)$ whenever $g\in
\ell_1(\Z,\mathcal{L}(E))$ and $u\in c_0(\Z,E)$.
\footnote{This follows immediately by approximating $g$ with the sequence 
$g_n:= \mathbb{1}_{[-n,+n]}g$, for $g_n\rightarrow g$ in $\ell_1$, 
$g_n *u\in c_0(\Z,E)$  and by Young's inequality $\|g*u-g_n*u\|_\infty =
\|(g-g_n)*u\|_\infty\leq \|g_n-g\|_1\|u\|_\infty \rightarrow 0$, so 
$g*u=\lim_{n\rightarrow \infty}g_n *u\in  c_0(\Z,E)$.}

Notice that the convolution with $g\in
\ell_1(\Z,\mathcal{L}(E))$ defines a bounded
linear operator $R_g$ on $c_0(E)$ by $u\mapsto g*u$ (more precisely, 
$R_g u = \rho(g*j(u))$).

\paragraph{Manifolds of sequences.} Let $M$ be a $C^k$ manifold modeled
on the Banach space $E$, let $x\in M$, and let $c_x(M)$ be the set
of sequences $u:\N \rightarrow M$ which converge to $x$. Equivalently,
denoting by $\overline{\N}=\N\cup \{\infty\}$ the one-point
compactification of the set of natural numbers, $c_x (M)$ is the set
\[
c_x(M) = \set{u\in C^0(\overline{\N},M)}{u(\infty)=x},
\]
so it can be endowed with the restriction of the compact-open 
topology of $C^0(\overline{\N},M)$. The space
$c_x(M)$ has the structure of a $C^k$ manifold modeled on the Banach
space $c_0(E)$. Indeed, given $C^k$ local charts $\varphi_n:U_n
\rightarrow \varphi_n(U_n)\subset E$, $n=0,\dots,m$, where $x\in U_m$
and $\varphi_m(x)=0$, consider the open subset of $c_x(M)$
\[
\mathcal{U} = 
\mathcal{U}(U_0,\dots,U_m) = \set{u\in c_x(M)}{u(n)\in U_n \; \forall
n=0,\dots,m-1, \; u(n)\in U_m \; \forall n\geq m},
\]
and the homeomorphism $\Phi=\Phi(\varphi_0,\dots,\varphi_m) : 
\mathcal{U} \rightarrow \Phi(\mathcal{U}) \subset c_0(E)$
defined by 
\[
\Phi(u)(n) = 
\varphi_n(u(n)) \mbox{ if } 0\leq n\leq m-1,
\quad \Phi(u)(n) = \varphi_m(u(n)) \mbox{ if } n\geq m. 
\]
It is easy to check that the collection of homeomorphisms
$\Phi(\varphi_0,\dots,\varphi_m)$ constitute a $C^k$ atlas of $c_x(M)$.
The tangent bundle of $c_x(M)$ is 
\[
Tc_x(M) = c_{0_x} (TM),
\]
where $0_x$ is the zero element of $T_x M\subset TM$, and its fibers
are
\[
T_u c_x (M) = \set{ v:\N \rightarrow TM}{v(n) \in T_{u(n)} M \;
  \forall n\in \N, \; \lim_{n\rightarrow \infty} v(n) = 0_x}.
\]
In particular, the tangent space of $c_x(M)$ at the constant sequence
$x$ is $T_x c_x(M) = c_0(T_x M)$.

The (left) shift operator $\mathcal{S}:c_x(M) \rightarrow c_x(M)$,
$\mathcal{S} (u)(n) = u(n+1)$, is of class $C^k$, and its differential
at $x$ is the (left) shift linear operator $S$ on $c_0(T_x M)$.

Also the evaluation at zero, $\mathrm{ev}_0:c_x(M) \rightarrow M$,
$u\mapsto u(0)$, is a map of class $C^k$, and its differential at $u$
is the linear evaluation at zero $D\mathrm{ev}_0(u)[v]=v(0)$. 

Finally, every continuous map $f:M \rightarrow N$ with 
$f(x)=y$ induces by composition a continuous map $f_*:c_x(M) 
\rightarrow c_y(N)$, $f_*(u):=f\circ u$.
If $f$ is of class $C^k$, $1\leq k \leq \omega$, 
so is $f_*$. Indeed, by local charts we are reduced to the case $M=E$,
$N=F$, $x=y=0$, where the $h$-th differential of $f_*$ at $u\in
c_0(E)$ is given by
the formula 
\[
\left(D^h f_*(u)[v]^h\right)(n)=D^h f(u(n)) [v(n)]^h.
\]
In particular, the differential of $f_*$ at the constant sequence $x$
is the multiplication operator by $Df(x)$,
\[
Df_*(x) : T_x c_x(M) = c_0(T_x M) \rightarrow T_y c_y (N) = c_0(T_y
N), \quad Df_*(x)[u](n) = Df(x)[u(n)].
\]

\paragraph{Hyperbolic fixed points.}
An invertible operator $T\in
\mathcal{L}(E)$ is said to be {\em hyperbolic} if its spectrum does not meet
the unit circle: $\sigma(T) \cap \{|z|=1\} = \emptyset$. Then
$\sigma(T)$ consists of the two disjoint closed subsets $\sigma(T)
\cap \{|z|<1\}$ and $\sigma(T)\cap \{|z|>1\}$, so $E$ has the
$T$-invariant spectral decomposition $E=E^s \oplus E^u$, 
where $\sigma(T|_{E^s}) =\sigma(T)
\cap \{|z|<1\}$ and  $\sigma(T|_{E^u}) =\sigma(T)\cap
\{|z|>1\}$. 

A fixed point $x$ of a diffeomorphism $f:M\rightarrow f(M)\subset 
M$ is said to be {\em
  hyperbolic} if the differential of $f$ at $x$, $Df(x)\in
\mathcal{L}(T_x M)$, is a hyperbolic operator. The corresponding
spectral decomposition of the tangent space at $x$ is denoted
by $T_x M = E^s \oplus E^u$.

\section{Proof of the stable manifold theorem}
 
Let us prove Theorem \ref{stab}.
By definition, $\mathcal{W}$ is the set
\[
\mathcal{W} = \set{u\in c_x(M)}{\mathcal{S}(u)=f_*(u)}.
\]
We start by studying the linear map $D\mathcal{S}(x) - Df_*(x) \in
\mathcal{L}(T_x c_x(M))$. By the discussion of section \ref{uno}, this
is the linear operator
\[
S-Df(x)_*: c_0(T_x M) \rightarrow c_0(T_x M).
\] 
Let us simplify the notation by setting $E=T_x M$, 
$T=Df(x)\in \mathcal{L}(E)$. Denote by $P^s$ and $P^u$ the spectral
projections associated to the decomposition $E=E^s\oplus E^u$.
The following lemma uses the fact that $T$ is a hyperbolic operator. 

\begin{lem}
\label{1}
Set, for $n\in\Z$,
\[
g(n):=T^{n-1}\left(\mathbb{1}_{\Z^+}(n) I_{E} -P^u \right),
\]
where $\Z^+=\{1,2,\dots\}$. Then $g\in
\ell_1(\Z,\mathcal{L}(E))$ and the corresponding convolution operator 
$R_g\in \mathcal{L}(c_0(E))$ is a right inverse of $S-T_*$. Moreover,
\[
\ker(S-T_*)=\set{u\in c_0(E)}{u(n)=T^n u(0) \; \forall 
n\in \N, \; u(0)\in E^s}.
\]
\end{lem}

\proof
Let $\|\cdot\|$ be the operator norm induced by a Banach norm on $E$. 
By the spectral radius theorem
\begin{eqnarray*}
\lim_{n\rightarrow \infty} \|g(n)\|^{1/n} = \lim_{n\rightarrow
  \infty} \|T^{n-1}P^s\|^{1/n} = \lim_{n\rightarrow \infty}
  \|T|_{E^s}^{n-1}\|^{1/n} = \max |\sigma(T|_{E^s})| < 1, \\
\lim_{n\rightarrow \infty} \|g(-n)\|^{1/n} = \lim_{n\rightarrow
  \infty} \|T^{-(n+1)}P^u\|^{1/n} = \lim_{n\rightarrow \infty}
  \|T|_{E^u}^{-(n+1)}\|^{1/n} = \max |\sigma(T|_{E^u}^{-1})| < 1.
\end{eqnarray*}
Therefore, $g(n)$ tends to 0 exponentially fast for $|n|\rightarrow
\infty$, in particular $g$ is in $\ell_1(\Z,\mathcal{L}(E))$. 

We have, for any $u\in c_0(E)$ and $n\in\N$,
\begin{eqnarray*}
\left[(S-T_*)R_g(u)\right](n)=\sum_{h=0}^{\infty}g(n+1-h)u(h)-
\sum_{h=0}^{\infty}Tg(n-h)u(h) \\
= \sum_{h=0}^{\infty}T^{n-h}\left[\mathbb{1}_{\Z^+}(n+1-h)-
\mathbb{1}_{\Z^+}(n-h)\right]u(h)=
\sum_{h=0}^{\infty}T^{n-h}\mathbb{1}_{\{0\}}(n-h)u(h)=u(n),
\end{eqnarray*}
that is, $(S-T_*)R_g=I_{c_0(E)}$.
\footnote{Here is a more heuristic argument to find a right inverse to
the linear operator $S-T_*$. First notice that the equation
$(S-T_*)u=w\in c_0(E)$ is equivalent to
$u(n+1)=Tu(n)+w(n)$, $\forall n\geq0$, that iterated gives
$u(n)=T^n u(0)+\sum_{h=0}^{n-1}T^{n-1-h}w(h)$. We can split this
equation into
$u(n)=T^n P^s u(0)+\sum_{h=0}^{n-1}T^{n-1-h}P^s w(h)+
T^n\left[ P^u u(0)+\sum_{h=0}^{n-1}T^{-1-h}P^u w(h)\right]$. 
Now the first and the second term converge as $n\rightarrow\infty$,
because the spectral radius theorem implies that 
$\|T^n P^s\|\leq c \lambda^n$, for some $c\geq 1$ and $\lambda<1$. 
The third term may not converge unless the
sequence into square brackets converges to $0$, that is, 
$P^u u(0)+\sum_{h=0}^{n-1}T^{-1-h}P^u w(h)=
-\sum_{h=n}^{\infty}T^{-1-h}P^u w(h)$, whence  
$u(n)=T^n P^s u(0)+(g*w)(n)$.}
Finally, it is clear that $u\in \ker(S-T_*)$ if and only if
$u(n)=T^n u(0)$ for any $n\in \N$, which defines an element of $c_0(E)$ 
if and only if $u(0)\in E^s$.
\qed

Let us prove that the closed subset $\mathcal{W}$ is a $C^k$
submanifold of $c_x(M)$. By Lemma \ref{1}, $D\mathcal{S}(x) - Df_*(x)$
is a linear retraction. Since the the space of linear retractions is
open, $D\mathcal{S}(u) - Df_*(u)$ is a linear retraction for every
$u\in \mathcal{W} \cap U$, for a suitable neighborhood $U$ of $x$ in
$c_x(M)$. 

On the other hand, $f_*$ and $\mathcal{S}$ commute. As a consequence if $u\in
\mathcal{W}$ then $u_n=f_*^n(u)$ is also in $\mathcal{W}$, and the linear
operator $D\mathcal{S}(u) - Df_*(u)$ is related to 
$D\mathcal{S}(u_n) - Df_*(u_n)$ by left and right multiplication by
invertible linear operators.
Since $u_n$ eventually belongs to $\mathcal{W}\cap U$, $D\mathcal{S}(u_n)
- Df_*(u_n)$ is a linear retraction, and so is
$D\mathcal{S}(u) - Df_*(u)$. Therefore, Proposition \ref{p1}
implies that $\mathcal{W}$ is a $C^k$ submanifold of $c_x(M)$, with
\begin{equation}
\label{tg}
T_x \mathcal{W} = \ker (D\mathcal{S}(x) - Df_*(x)) = \set{v\in c_0(T_x
  M)}{v(n) = Df(x)^n v(0)\; \forall n\in \N, \; v(0)\in E^s}.
\end{equation}

The $C^k$ map $\mathrm{ev}_0:c_x(M) \rightarrow M$ is clearly
injective on $\mathcal{W}$, and it remains to show that for any
$u\in \mathcal{W}$, $D \mathrm{ev}_0|_{\mathcal{W}} (u):T_u
\mathcal{W} \rightarrow T_{u(0)} M$ is a linear section.
This is clearly true if $u=x$, for the expression (\ref{tg}) shows that
$D\mathrm{ev}_0(x)$ is an isomorphism onto $E^s$, which splits in
$T_x M$. Since the space of linear sections is open, the same is true
for every $u$ in a neighborhood of $x$ in $\mathcal{W}$. 
The formula $\mathrm{ev}_0 \circ f_* = f \circ \mathrm{ev}_0$ implies
that $D\mathrm{ev}_0(u)|_{T_u \mathcal{W}}$ is obtained from
$D\mathrm{ev}_0 (u_n)|_{T_{u_n} \mathcal{W}}$ by left and
right multiplication by invertible operators. As before, since $u_n$
converges to $x$, we conclude that  $D\mathrm{ev}_0(u)|_{T_u
  \mathcal{W}}$ is a linear section for every $u\in \mathcal{W}$. 
The proof of Theorem \ref{stab} is complete.
 
\begin{rem}
The fact that $D\mathcal{S}(u)-Df_*(u)$ is a linear retraction can
also be proved directly by the following generalization of Lemma
\ref{1}: if $\mathcal{T}:c_0(E)\rightarrow c_0(E)$ is the
multiplication operator by a sequence $(T_n)\subset \mathcal{L}(E)$
converging to a hyperbolic operator, then $S-\mathcal{T}\in
\mathcal{L}(c_0(E))$ is a linear retraction.
\end{rem} 

\begin{rem}
\label{rof}
A similar argument allows to prove the stable manifold theorem for 
a hyperbolic equilibrium point $x$ of the local flow  
determined by a $C^k$ vector field $X$ on a $C^{k+1}$ Banach manifold
$M$, where $1\leq k\leq \omega$. 
Denote by $C^1_x([0,+\infty[,M)$ the space
of $C^1$ curves $[0,+\infty[\rightarrow M$ converging to $x$ for
$t\rightarrow +\infty$, with the first derivative converging
to 0. Then one can use the implicit function theorem to prove that the set
\[
\mathcal{W} = \set{u\in C^1_x([0,+\infty[,M)}{u^{\prime}(t)-X(u(t))=0}
\]
is a $C^k$ submanifold of $C^1_x([0,+\infty[,M)$, and that the
evaluation at $0$ subordinates a $C^k$ immersion of $\mathcal{W}$ onto
the stable manifold of $0$. The basic linear step consists in proving that
if $L\in \mathcal{L}(E)$ is infinitesimally hyperbolic (i.e.\ its
spectrum does not meet the imaginary axis), then the operator
\[
\frac{d}{dt} - L : C_0^1([0,+\infty[,E) \rightarrow
C_0^0([0,+\infty[,E)
\]
is a linear retraction. See \cite{ama04m} for more details. 
\end{rem}

\section{The local stable manifold theorem}
\label{tre}

\paragraph{Adapted norms.}
An {\em adapted norm} for a hyperbolic operator $T\in \mathcal{L}(E)$ 
is an equivalent norm $|\cdot|$ on $E$ such that 
\begin{equation}
\label{ada}
|\xi| = \max\{ |P^s \xi|, |P^u \xi|\}, \quad
|T\xi | \leq \lambda |\xi| \quad \forall \xi \in E^s, \quad  
|T^{-1}\xi | \leq \lambda |\xi| \quad \forall \xi \in E^u,
\end{equation}
for some $0<\lambda<1$. One can always find a norm with these
properties, provided that $\lambda> \max | \sigma(T) \cap \{|z|<1\}|$
and $\lambda > 1/\min |\sigma(T)\cap \{|z|>1\}|$. Indeed, the
following stronger statement holds: if the spectrum of $T$ is
contained in the annulus $\{\alpha < |z| < \beta\}$, then $E$ has an
equivalent norm $|\cdot|$ such that, denoting by $\|\cdot\|$ the corresponding
operator norm, there holds
\begin{equation}
\label{adno}
\|T\| \leq \beta, \quad \|T^{-1}\| \leq \frac{1}{\alpha}.
\end{equation}
If $|\cdot|_0$ is any equivalent norm on $E$, a norm $|\cdot|$
satisfying (\ref{adno}) can be defined as 
\[
|\xi| = \sum_{n=1}^{\infty} \alpha^{n} |T^{-n} \xi|_0 +
 \sum_{n=0}^{\infty} \beta^{-n} |T^n \xi|_0,
\]
as shown by the the spectral radius theorem (see for instance
\cite{hps77}, Proposition 2.8). 

\paragraph{The local stable manifold.} Let $U$ be an open
neighborhood of 0 in the Banach space $E$, and let $f:U\rightarrow
f(U)\subset E$ be a $C^k$ diffeomorphism, $1\leq k \leq \infty$ or
$k=\omega$, having $0$ as a hyperbolic fixed point. Let $T=Df(0)$, let
$E=E^s \oplus E^u$ be the corresponding splitting, and let $|\cdot|$
be an adapted norm on $E$ for the hyperbolic operator $T$. If $V$ is a
closed linear subspace of $E$, we denote by $V(r)$ the closed
ball in $V$ of radius $r$. By the first property of adapted norms (\ref{ada}),
$E(r)=E^s(r) \times E^u(r)$.

Given $r>0$ such that $E(r)\subset U$, the {\em local stable manifold}
of $0$ is the set
\[
W^s_{\mathrm{loc},r} (0) := \set{p\in \bigcap_{n\in \N} f^{-n}
  (E(r))}{\lim_{n\rightarrow \infty} f^n(p)=0}.
\]
This definition depends on $r$. However, if $r_0$ is small enough
\begin{equation}
\label{varr}
W^s_{\mathrm{loc},r} (0) = W^s_{\mathrm{loc},r_0}(0) 
\cap E(r) \quad \forall r \leq r_0.
\end{equation}
Indeed, the first set is contained in the second one by definition. 
Let us prove that the other inclusion holds when $r_0$ is small. 
The point $0\in \mathcal{W}$ is a fixed point of the $C^k$ map
$f_*|_{\mathcal{W}} : \mathcal{W} \rightarrow \mathcal{W}$.
By (\ref{tg}) and by the second property of adapted norms (\ref{ada}),
\[
\|Df_*(0)|_{T_0 \mathcal{W}}\| = 
\|T_*|_{T_0 \mathcal{W}} \|\leq \lambda < 1,
\]
so a first order approximation shows that $f_*|_{\mathcal{W}}$ is locally a 
contraction at $0$. In particular, there exists $r_0>0$ such that 
$\|f_*(u)\|_{\infty}
  \leq \|u\|_{\infty}$ for every $u\in \mathcal{W}$ with
  $\|u\|_{\infty} \leq r_0$. Since $f_*$ coincides with the shift $S$ on
  $\mathcal{W}$, this is equivalent to saying that $|u(n)| =
  |f^n(u(0))|$ is a decreasing sequence, if $u\in \mathcal{W}$ and
  $\|u\|_{\infty} \leq r_0$. The conclusion follows.

\begin{thm} {\em (Local stable manifold theorem)} Let $U$ be an open
neighborhood of 0 in the Banach space $E$, $f:U \rightarrow
f(U)\subset E$ be a $C^k$ diffeomorphism, $1\leq k \leq \omega$, 
having $0$ as a hyperbolic fixed point, inducing the
splitting $E=E^s \oplus E^u$. 
If $r>0$ is small
enough, then $W^s_{\mathrm{loc},r}(0)$ is the graph of a $C^k$ map 
$w: E^s(r) \rightarrow E^u(r)$ such that $w(0)=0$ and $Dw(0)=0$.
\end{thm}  

\proof
Since $\mathrm{ev}_0|_{\mathcal{W}}$ is a $C^k$ immersion, it is an
embedding of an open neighborhood of $0$ in $\mathcal{W}$. Therefore, if
$r_0$ is small enough the set
\[
\set{\mathrm{ev}_0(u)}{u\in \mathcal{W},
  \; \|u\|_{\infty} < r_0}  
\]
is a $C^k$ submanifold of $E$, with tangent space at $0$
$D\mathrm{ev}_0(0) T_0 \mathcal{W} =
E^s$. Then, if $r< r_0$ is small enough, the set 
\[
W^s_{\mathrm{loc},r_0}(0) \cap E(r) =
\set{\mathrm{ev}_0(u)}{u\in \mathcal{W},
  \; \|u\|_{\infty} < r_0} \cap E(r)
\]
is the graph of a $C^k$ map $w:E^s(r) \rightarrow
E^u(r)$ such that $w(0)=0$ and $Dw(0)=0$, and the conclusion follows
from (\ref{varr}).
\qed  

\section{The diffeomorphism class of the global stable manifold}

\paragraph{Characterization of diffeomorphisms such that $\mathbf{W^s(x)}$ is
  an immersed copy of $\mathbf{E^s}$.}
A hyperbolic fixed point $x$ of a diffeomorphism $f:M\subset
M^{\prime} \rightarrow f(M)\subset M^{\prime}$ is said to be a {\em local
  attractor} if $W^s(x)$ is a neighborhood of $x$, or
equivalently\footnote{Indeed, if $E^s \neq T_x M$ the local stable
manifold has empty interior, so by Baire theorem the global stable
manifold cannot fill an open set.} if $E^s=T_x M$. 
It is said to be a {\em global attractor} if $W^s(x)=M$ (in
particular, $f(M)\subset M$), in which case $f$ is said to be a {\em
topological contraction of $M$}.  

\begin{defn}
\label{prima}
Let $U$ be an open neighborhood of $0$ in the Banach space $E$, and
let $f:U \rightarrow f(U)\subset E$ be a diffeomorphism of class
$C^k$, $1 \leq k \leq \omega$, such that the hyperbolic
fixed point $0$ is a local attractor. 
We say that the germ of $f$ at $0$ extends to a $C^k$ topological
contraction of $E$ if it can be represented by a global $C^k$ 
diffeomorphism of $E$ having $0$ as a global attractor. 
\end{defn}

In other words, we are asking the existence of a global
diffeomorphism $\tilde{f} :E \rightarrow E$ of class $C^k$ which
coincides with $f$ in a neighborhood $V\subset U$ of $0$, and which
is a topological contraction. If we ask $\tilde{f}$ to be only
Lipschitz, its existence is always guaranteed:

\begin{lem}
\label{lippo}
Under the assumptions of Definition \ref{prima}, there always 
exists a homeomorphism $\tilde{f}: E \rightarrow E$ such that
$\tilde{f}$ and $\tilde{f}^{-1}$ are globally Lipschitz, $\tilde{f}$
coincides with $f$ in a neighborhood $V\subset U$ of $0$, and such that
$0$ is a global attractor. 
\end{lem}

\proof
Let $T=Df(0)$ and let $|\cdot|$ be an adapted norm 
for $T$: if $\|\cdot\|$ denotes the corresponding operator norm, we have
$\|T\|<1$. Let $\varphi: E \rightarrow E$
be a $k$-Lipschitz bounded map which coincides with the identity in a
neighborhood of $0$ (for instance, $\varphi(\xi)=\chi(\xi) \xi$, with
$\chi(s)=1$ for $s\leq 1$, $\chi(s)=2-s$ for $1\leq s \leq 2$,
$\chi(s)=0$ for $s\geq 2$, has the required properties, with $\lip
\varphi \leq 3$).  Write $f$ as $f(\xi) = T(\xi + f_0(\xi))$. Since $f_0$ is
at least $C^1$ and $Df_0(0)=0$, we can find a neighborhood $U_0$ of
$0$ such that
\[
\lip f_0|_{U_0} < 1/k, \quad \theta:= (1 + k\, \lip f_0|_{U_0}) \|T\|<1.
\]
Up to replacing $\varphi(\xi)$ by $\lambda \varphi(\xi/\lambda)$ -
which is still $k$-Lipschitz - 
we may assume that $\varphi(E)\subset U_0$, and we set 
\[
\tilde{f} : E \rightarrow E, \quad \tilde{f}(\xi) = T (\xi + f_0
(\varphi(\xi))).
\]
The map $\tilde{f}$ is a Lipschitz diffeomorphism of $E$ onto $E$
together with its inverse
because $\lip (f_0 \circ \varphi) \leq k\, \lip f_0|_{U_0}<1$. Moreover,
$\tilde{f}=f$ in a neighborhood of $0$, and
\[
|\tilde{f}(\xi)| \leq \|T\| (|\xi| + |f_0(\varphi(\xi))|) \leq \|T\| (
1 + k\, \lip f_0|_{U_0} ) |\xi| = \theta |\xi|,
\]
so the fact that $\theta<1$ implies that $0$ is a global attractor.
\qed

In the analytic category the extension of a germ is unique (and not
always possible, even in finite dimensional spaces), so the request 
of Definition \ref{prima} is quite strong. A weaker property is given
by the following:

\begin{defn}
Under the same assumptions of Definition \ref{prima}, we say that
the germ of $f$ at $0$ extends up to conjugacy 
to a $C^k$ topological contraction of $E$ if there exists a 
$C^k$ diffeomorphism 
$h:V\rightarrow h(V)$ of an open neighborhood of $0$ with $h(0)=0$ 
such that the germ of $h\circ f \circ h^{-1}$ at $0$ extends to a
topological contraction of $E$.
\end{defn} 

For instance, if an analytic diffeomorphism is analytically
linearizable near a hyperbolic local attractor, then it extends up to
conjugacy to an
analytic topological contraction of $E$. This latter definition is
relevant also in the smooth and in the finite differentiability
category because it extends to diffeomorphisms defined on manifolds:

\begin{defn}
If $U$ is an open neighborhood of $x$ in the Banach manifold $M$ 
 - modeled on the Banach space $E$ - and
$f:U \rightarrow f(U)\subset M$ is a $C^k$ diffeomorphism, $1\leq
 k\leq \omega$, having $x$ as a local attractor, we
say that the germ of $f$ at $x$ extends up to conjugacy 
to a $C^k$ topological
contraction of $E$ if the germ of diffeomorphism 
on an open neighborhood of $0$ in $E$ defined by conjugacy by a local chart 
(mapping $x$ into $0$) extends up to conjugacy 
to a topological contraction of $E$.
\end{defn} 

Clearly, this definition does not depend on the choice of the local
chart. 

\begin{thm}   
\label{gs}
Let $f$ be a $C^k$ diffeomorphism of the $C^k$ Banach manifold $M$
modeled on the Banach space $E$, $1\leq k \leq \omega$, 
and assume that the hyperbolic fixed point $x\in M$ is a global attractor.
Then $M$ is homeomorphic to $E$ by a bi-locally Lipschitz homeomorphism.
Furthermore, $M$ is $C^k$ diffeomorphic to $E$ if and only if the germ
of $f$ at $x$ extends up to conjugacy to a $C^k$ topological
contraction of $E$.
\end{thm}

\proof
Let $U\subset M$ be an open neighborhood of $x$ in $M$, and let
$\psi: U \rightarrow \psi(U)\subset E$ be a $C^k$ local chart
with $\psi(x)=0$. Consider the $C^k$ diffeomorphism
\[
g = \psi \circ f \circ \psi^{-1} : \psi(U\cap f^{-1}(U))
\rightarrow \psi(U \cap f(U))\subset E.
\]
By Lemma \ref{lippo}, there exists an invertible topological contraction 
$h:E \rightarrow E$ which is locally Lipschitz together with its 
inverse and which coincides with $g$ in an open neighborhood
$V\subset  \psi(U\cap f^{-1}(U))$ of $0$.
If we are also assuming that the germ of $f$ at $x$ extends up to
conjugacy to a
$C^k$ topological contraction of $E$, up to changing the 
local chart $\psi$ we may assume that $h$ is a $C^k$ diffeomorphism.

Let us define the global homeomorphism $\phi:M \rightarrow E$. Let
$p\in M$. Since $W^s(x)=M$, there exists $n\in \N$ such that
$f^n(p)\in \psi^{-1}(V)$, and we set
\[
\phi(p) = h^{-n} ( \psi (f^n(p))).
\]
Since $\psi$ conjugates the diffeomorphisms $f|_{\psi^{-1}(V)}$
and $h|_V = g|_V$, this definition does not depend on the choice of
$n$. The map $\phi$ is invertible, its inverse
being the map
\[
\phi^{-1} (\xi) = f^{-n} (\psi^{-1} (h^n(\xi))),
\]
for $n=n(p)$ so large that $h^n(\xi) \in V$. So $\phi$ is the required
locally Lipschitz homeomorphism, and when $h$ is a $C^k$
diffeomorphism, $\phi$ is also a $C^k$ diffeomorphism.

Conversely, assume that there is a $C^k$ diffeomorphism $\phi: M \rightarrow
E$. Up to composition with a translation, we may assume that
$\phi(x)=0$. In particular, $\phi$ is a local chart mapping $x$
into $0$, and $\phi\circ f \circ \phi^{-1}$ is a global
diffeomorphism of $E$ onto $E$ which is a topological contraction of
$E$. This shows that the germ of $f$ at $x$ extends up to conjugacy 
to a $C^k$ topological contraction of $E$. 
\qed   

Let us consider the general case $W^s(x)\neq M$. Denote by
$W^s_{\mathrm{loc}}(x)$ the image by $\mathrm{ev}_0$ of a neighborhood
of $x$ in $\mathcal{W}$, so small that $\mathrm{ev}_0$ is an embedding
on it. Then $W^s_{\mathrm{loc}}(x)\subset W^s(x)$ is a
$C^k$-submanifold modeled on the Banach space $E^s$. 

\begin{cor}
\label{dopo}
Let $f:M \rightarrow M$ be a $C^k$ diffeomorphism of a Banach
manifold, $1\leq k \leq \omega$. 
Let $x$ be a hyperbolic fixed point of $f$, with associated
splitting $T_x M = E^s \oplus E^u$. Then the $C^k$ manifold $\mathcal{W}$
is homeomorphic to $E^s$ by a bi-locally Lipschitz homeomorphism.
Furthermore, $\mathcal{W}$ is $C^k$ diffeomorphic to $E^s$ if and only 
if the germ of $f|_{W^s_{\mathrm{loc}}(x)}$ at $x$ extends up to conjugacy 
to a $C^k$ topological contraction of $E^s$. 
\end{cor}   

Indeed, it is enough to apply Theorem \ref{gs} to the 
diffeomorphism $f_*$ of the manifold $\mathcal{W}$.

Since the evaluation at zero defines a local conjugacy between the
restriction of $f_*$ to a small neighborhood of $x$ in $\mathcal{W}$
and the restriction of $f$ to $W^s_{\mathrm{loc}}(x)$, the above corollary
can be restated in the following way: $W^s(x)$ is always the image of
a locally Lipschitz and locally closed injective map $E^s\hookrightarrow M$, 
which can be chosen to be also a $C^k$ immersion if and only if the germ of
$f|_{W^s_{\mathrm{loc}}(x)}$ at $x$ extends up to conjugacy 
to a $C^k$ topological contraction of $E^s$.   

The statement about the locally Lipschitz homeomorphism in the above
corollary is the first part of Theorem \ref{diffcl}.

\paragraph{Regularity conditions on the Banach space $\mathbf{E}$.}
If $1\leq k \leq \infty$ and the Banach space is somehow regular, one 
may expect that every germ of $C^k$ diffeomorphism having $0$ as a 
hyperbolic local attractor extends to a $C^k$ topological contraction of $E$. 
This is actually true for a large class of Banach spaces. 
Indeed let us consider the following conditions on $E$:
\renewcommand{\theenumi}{\arabic{enumi}}
\renewcommand{\labelenumi}{(E\theenumi)}
\begin{enumerate}
\item $E$ is finite dimensional;
\item $E$ has a Hilbert structure;
\item there exists a $C^k$ norm on $E$ (i.e.\ an equivalent norm
  $|\cdot|$ such that the function $\xi \mapsto |\xi|$ is $C^k$ on
  $E\setminus \{0\}$);
\item there exists a bounded $C^k$ map $\varphi:E \rightarrow E$ which
  coincides with the identity on a neighborhood of $0$;
\item there exists a bounded $C^k$ map $\varphi:E \rightarrow E$ such
  that $\varphi(0)=0$ and $D\varphi(0)=I$.
\end{enumerate}
\renewcommand{\theenumi}{\roman{enumi}}
\renewcommand{\labelenumi}{(\theenumi)}
It is readily seen that (E1) $\implies$ (E2) $\implies$ (E3) $\implies$
(E4) $\implies$ (E5). For instance, if $E$ admits a $C^k$ norm 
$| \cdot |$, a $C^k$ map $\varphi$ satisfying (E4) can be defined as
\[
\varphi( \xi ) = \chi (|\xi|) \xi,
\]
where $\chi:[0,+\infty) \rightarrow \R$ is a smooth function with
compact support and such that $\chi(s)=s$ for $s\leq 1$. Condition
(E4) was considered by Atkin \cite{atk01}, who observed that such a
condition may hold also for Banach spaces which do not have regular
norms.
For instance, the Banach space $\ell_{\infty} =
\ell_{\infty}(\N,\R)$ does not admit a smooth norm, but it supports a
smooth map $\varphi$ satisfying (E4), namely
\[
\varphi (u) (n) = \chi(u(n)), \quad \forall n\in \N,
\]
where $\chi:\R \rightarrow \R$ is a smooth bounded function coinciding
with the identity in a neighborhood of $0$. A similar construction
works for the Banach space $L^{\infty}(X,\mathcal{F},\mu)$,
$(X,\mathcal{F},\mu)$ a measure space, and for the Banach space
$C^0(K,\R)$, $K$ a compact topological space.

The condition (E5) holds for a large class of Banach
spaces, even in the case $k=\omega$, in which (E4) is obviously never
fulfilled. For instance, the spaces $L^{\infty}(X,\mathcal{F},\mu)$ and
$C^0(K,\R)$ admit the analytic map
\[
u \mapsto \sin u,
\]
which satisfies (E5). It is not clear whether (E4) and (E5) are
equivalent for $k<\omega$.  

Even if we do not know any counterexample, it is not likely that
every Banach space admits a smooth map satisfying (E4) or (E5). 
In any case, it would be interesting to
characterize the class of Banach spaces which admit such maps. At the
moment the situation is unclear even for simple spaces such as $\ell_1$.  

\paragraph{The case of finite differentiability.}
The following result, together with Corollary \ref{dopo}, implies the
(i) part of Theorem \ref{diffcl}. 

\begin{prop}
\label{finreg}
Let $1\leq k \leq \infty$, and assume that the Banach space $E$
satisfies condition (E4) above. Let $U\subset E$ be an 
open neighborhood of $0$, and let  $f:U \rightarrow f(U)\subset E$ be a $C^k$
diffeomorphism with hyperbolic fixed point $0$ which is a local attractor.
Then the germ of $f$ at $0$ extends to a $C^k$ topological contraction
of $E$.
\end{prop}

Notice that if the map $\varphi$ appearing in (E4) 
is also globally Lipschitz, we could argue 
as in Lemma \ref{lippo}, writing $f =T\circ (\mathrm{id} + f_0)$ and then
using $\varphi$ to extend $f_0$ to a $C^k$ map on $E$ with small
Lipschitz norm. Without this assumption, a natural idea is to see
$\mathrm{id} + f_0$ as the time-one map obtained by integrating a time
dependent small vector field $X$, and then use $\varphi$ to extend $X$. We
need therefore the following easy:

\begin{lem}
\label{campi}
Let $1\leq k \leq \omega$, let $U\subset E$ be an open neighborhood of
$0$, and let $g:U \rightarrow E$ be a $C^k$ map such that $g(0)=0$ and
$Dg(0)=I$. Then there exists a neighborhood $V\subset U$ of $0$ and a
$C^k$ map $X:[0,1]\times V \rightarrow E$ such that $X(t,0)=0$,
$D_2X(t,0)=0$ for every $t\in [0,1]$, and such that the solution of the
Cauchy problem
\[
\partial_t G(t,\xi) = X(t,G(t,\xi)), \quad G(0,\xi) = \xi,
\]
satisfies $G(1,\xi)=g(\xi)$ for every $\xi$ in some neighborhood of
$0$.
\end{lem}  

\proof
The differential with respect to the second variable of the $C^k$ map
\[
G : [0,1] \times U \rightarrow E, \quad G(t,\xi) = tg(\xi) + (1-t)
\xi,
\]
namely $t Dg(\xi) + (1-t) I$, is uniformly invertible for every $(t,\xi)$ in a
neighborhood of $[0,1]\times \{0\}$. By the parametric inverse 
mapping theorem, there exist a neighborhood $V$ of $0$ and a $C^k$ map 
$H: [0,1] \times V  \rightarrow E$ such that
\[
H(t,G(t,\xi)) = \xi \quad \forall (t,\xi) \in [0,1]\times V.
\]
Setting
\[
X(t,\eta) = g(H(t,\eta)) - H(t,\eta),
\]
we get that $X$ has the desired properties.
\qed

We can now prove Proposition \ref{finreg}.

\proof[of Proposition \ref{finreg}]
Let $T=Df(0)$ and let $|\cdot|$ be an adapted norm for $T$, so that
$T$ becomes a contraction. Consider the diffeomorphism $g=T^{-1} f$, whose
differential at $0$ is the identity operator. Let $X\in C^k( [0,1]
\times V,E)$ and $G$ be as in Lemma \ref{campi}. 
By assumption, there is a $C^k$ map $\varphi:E \rightarrow E$ whose
image is contained in $B_{r_0}(0)$, which coincides with the identity
map on $B_{s_0}(0)$, for some $0<s_0<r_0<+\infty$. Since $X(t,0) = 0$
and $D_2 X(t,0)=0$, for every $t\in [0,1]$, we can find $r>0$ such
that
\begin{equation}
\label{linest}
|X(t,\xi)| \leq \epsilon |\xi|, \quad \forall (t,\xi) \in [0,1] \times
B_r(0),
\end{equation}
where $\epsilon>0$ is so small that $e^{\epsilon r_0/s_0} \|T\|<1$.
The $C^k$ map 
\[
\psi(\xi) = \frac{r}{r_0} \varphi\left( \frac{r_0}{r} \xi \right),
\]
takes values in $B_r(0)$, and coincides with the identity mapping on
$B_s(0)$, with $s=rs_0/r_0$. The time-dependent $C^k$ vector field
\[
\tilde{X}: [0,1]\times E \rightarrow E, \quad \tilde{X}(t,\xi ) 
= X(t,\psi(\xi)),
\]
coincides with $X$ on $[0,1]\times B_s(0)$, so by (\ref{linest}) 
$|\tilde{X}(t,\xi)| \leq \epsilon |\xi|$ there. On the other hand, if
$|\xi|\geq s$,
\begin{equation}
\label{cresc}
|\tilde{X}(t,\xi)| = |X(t,\psi(\xi))| \leq \epsilon |\psi(\xi)| \leq
\epsilon r \leq \epsilon \frac{r}{s} |\xi| = \epsilon \frac{r_0}{s_0}
|\xi|.
\end{equation}
We conclude that the above estimate holds for every $(t,\xi)\in
[0,1]\times E$. Therefore, the solution $\tilde{G}$ of the Cauchy problem
\[
\partial_t \tilde{G}(t,\xi) = \tilde{X}(t,\tilde{G}(t,\xi)), 
\quad \tilde{G}(0,\xi) = \xi,
\]
is defined for every $(t,\xi)\in [0,1]\times E$, coincides with
$G$ in a neighborhood of $[0,1]\times \{0\}$, and by (\ref{cresc}) it 
satisfies
\begin{equation}
\label{ultima}
|\tilde{G}(t,\xi)| \leq e^{\epsilon \frac{r_0}{s_0} t} |\xi| , \quad
\forall (t,\xi) \in [0,1]\times E.
\end{equation}
Since $G(1,\xi)=g(\xi)$, the global $C^k$ diffeomorphism
\[
\tilde{g}: E \rightarrow E, \quad \tilde{g}(\xi) = \tilde{G}(1,\xi),
\]
coincides with $g$ in a neighborhood of $0$. Then the $C^k$
diffeomorphism $\tilde{f} = T \tilde{g}$ coincides with $f$ in a
neighborhood of $0$. Finally, by (\ref{ultima}),
\[
|\tilde{f}(\xi)| \leq \|T\| \, |\tilde{g}(\xi)| \leq \|T\| e^{\epsilon
  r_0/s_0} |\xi|,
\]
and the fact that $e^{\epsilon r_0/s_0} \|T\|<1$ implies that the
hyperbolic fixed point $0$ is a global attractor. 
\qed   

\paragraph{The smooth and the analytic cases.}
It remains to prove the (ii) part of Theorem \ref{diffcl}. Let us
start by examining some consequences of assumption (E5) (here by Taylor
polynomial we mean Taylor polynomial at $0$):

\begin{lem}
\label{equi}
Let $k=\infty$ or $k=\omega$. The following facts are equivalent:
\begin{enumerate}
\item there exists a bounded map $\varphi\in C^k(E,E)$ such that
  $\varphi(0)=0$ and $D\varphi(0)=I$;
\item every polynomial map $p:E\rightarrow E$, $\deg p\leq n$, is the
  Taylor polynomial of order $n$ of some bounded map $\psi\in C^k
  (E,E)$;
\item if $\epsilon>0$, $K$ is a compact manifold of class $C^h$, 
  $0\leq h \leq
  \omega$, $U$ is an open neighborhood of $0$ in $E$, $F$ is a Banach 
  space, and $\phi: K \times U \rightarrow F$
  is a $C^h$ map such that for every $x\in K$ the map $\phi(x,\cdot) : E
  \rightarrow F$ is of class $C^k$ and $\phi(x,0)=0$, then there exists a map
  $\psi: K \times E \rightarrow F$ with the same regularity, such that
  $\sup |\psi|<\epsilon$, and such
  that for every $x\in K$ the Taylor polynomials of order $n$ 
  of  $\phi(x,\cdot)$ and of $\psi(x,\cdot)$ coincide.
\end{enumerate}
\end{lem}
   
\proof
Statement (i) is a particular case of (ii): take $p(x)=x$ and
$n=1$. Statement (ii) is a particular case of (iii): take $K$ a
singleton, $\phi=p-p(0)$. So it is enough to prove that (i) implies
(iii). Given $\varphi$ satisfying (i), it is easy to construct a
bounded map $\varphi_n\in C^k(E,E)$ such that $\varphi_n(\xi) = \xi +
o(|\xi|^n)$ for $\xi\rightarrow 0$. Indeed, one may define $\varphi_n$
inductively as
\[
\left\{ \begin{array}{l} \varphi_1 = \varphi, \\ \varphi_{n+1} (\xi) =
  \varphi_n \bigl( \xi - \frac{1}{(n+1)!} D^{n+1} \varphi_n (0)
  \xi^{n+1} \bigr). \end{array} \right.
\]
Let $\phi: K \times U \rightarrow F$ be the map appearing in (iii).
By replacing $\varphi(\xi)$ by $\lambda \varphi (\xi/\lambda)$ in the
above construction, we may assume that the image of $\varphi$ - hence
also the image of $\varphi_n$ - is contained in a
neighborhood $V\subset U$ of $0$ which is so small that $\sup_{K
  \times V} |\phi|<\epsilon$ (such a neighborhood exists because $\phi$ is
continuous, and it vanishes on the compact set $K\times \{0\}$). Then the map
\[
\psi(x,\xi) = \phi(x,\varphi_n(\xi))
\]
has the required properties.
\qed   

The proof of Proposition \ref{411} below relies on the following classical
conjugacy result (see for instance \cite{cfl03}):
\begin{thm} 
\label{conju}
Let $k=\infty$ or $k=\omega$. Let $U\subset E$ be an open
neighborhood of $0$, and let $f: U \rightarrow f(U)\subset E$ be a
$C^k$ map with $f(0)=0$, $Df(0)=T$ an isomorphism with 
spectral radius $\rho(T)<1$. Then $f$ is $C^k$ locally conjugated to
its Taylor polynomial of order $n$, provided that $n$ is so large that
$\rho(T^{-1}) \rho(T)^{n+1}<1$.
\end{thm} 

The (ii) part of Theorem \ref{diffcl} is a consequence of Corollary
\ref{dopo} and of the following:

\begin{prop}
\label{411}
Let $k=\infty$ or $k=\omega$, and assume that the Banach space $E$
satisfies condition (E5). Let $U\subset E$ be an open neighborhood of
$0$, and let $f:U \rightarrow f(U)\subset E$ be a $C^k$ diffeomorphism
with hyperbolic fixed point $0$, which is a local attractor. Then the
germ of $f$ at $0$ extends up to conjugacy 
to a $C^k$ topological contraction of $E$.
\end{prop}

\proof
Let $T=Df(0)$ and let $|\cdot|$ be an adapted norm for $T$, so that
$T$ becomes a contraction. Consider the diffeomorphism $g:= T^{-1} f$,
whose differential at $0$ is the identity operator. Let $X\in
C^k([0,1]\times V,E)$ and $G$ be as in Lemma \ref{campi}: $G$ is the
map obtained by integrating the time dependent vector field $X$, and
$G(1,\cdot)=g$ in a neighborhood of $0$. 

Consider the $C^k$ map 
\[
\phi: [0,1]\times V \rightarrow \mathcal{L}(E), \quad \phi(t,\xi) = D_2
X(t,\xi).
\]
Fix $n\in \N$ so large that 
\begin{equation}
\label{ro}
\rho(T^{-1}) \rho(T)^{n+1}<1,
\end{equation}
and $\epsilon>0$ so small that $e^{\epsilon} \|T\|<1$. Since
$\phi(t,0) = D_2X (t,0)=0$, by Lemma \ref{equi} there exists a $C^k$
map $\psi:[0,1] \times E \rightarrow L(E)$ such that
\begin{equation}
\label{limi}
\sup_{(t,\xi) \in [0,1]\times E} \|\psi(t,\xi)\| < \epsilon,
\end{equation}
and such that for every $t\in [0,1]$ the Taylor polynomials of order
$n-1$ of $\phi(t,\cdot)$ and $\psi(t,\cdot)$ coincide. Equivalently:
\begin{equation}
\label{coinc}
\psi(t,\xi) - \phi(t,\xi) = o(|\xi|^{n-1}) \quad \mbox{for } \xi
\rightarrow 0,
\end{equation}
uniformly in $t\in [0,1]$. Consider the globally defined time
dependent vector field of class $C^k$,
\[
\tilde{X} : [0,1]\times E \rightarrow E, \quad \tilde{X}(t,\xi) =
\int_0^1 \psi(t,s\xi)\xi\, ds.
\]
By (\ref{limi}),
\begin{equation}
\label{croce}
|\tilde{X}(t,\xi)| \leq \int_0^1 \|\psi(t,s\xi)\|\, |\xi|\, ds \leq
\epsilon |\xi|, \quad \forall (t,\xi) \in [0,1] \times E.
\end{equation}
Since 
\[
X(t,\xi) = \int_0^1 \frac{d}{ds} X(t,s\xi) \, ds = \int_0^1 D_2
X(t,s\xi) \xi \, ds = \int_0^1 \phi(t,s\xi) \xi \, ds,
\]
by (\ref{coinc}) there holds
\[
\tilde{X}(t,\xi) - X(t,\xi) = \int_0^1 (\psi(t,s\xi) - \phi(t,s\xi))
\xi\, ds = o(|\xi|^n) \quad \mbox{for } \xi
\rightarrow 0,
\]
uniformly in $t\in [0,1]$. Therefore, $X(t,\cdot)$ and
$\tilde{X}(t,\cdot)$ have the same Taylor polynomial of order $n$. An
easy induction argument then shows that $G(t,\cdot)$ and the solution
$\tilde{G}(t,\cdot)$ of the Cauchy problem
\[
\partial_t \tilde{G}(t,\xi) = \tilde{X}(t,\tilde{G}(t,\xi)), \quad
\tilde{G}(0,\xi)=0, 
\]
have the same Taylor polynomial of order $n$, for every $t\in
[0,1]$. In particular, $g(\xi)=G(1,\xi)$ and $\tilde{g}(\xi) :=
\tilde{G}(1,\xi)$ have the same Taylor polynomial of order $n$. The
same happens for $f=T \circ g$ and $\tilde{f} := T \circ \tilde{g}$,
so by (\ref{ro}) Theorem \ref{conju} implies that $f$ and $\tilde{f}$
are $C^k$ locally conjugated at $0$.   

Since $\tilde{X}$ has linear growth by
(\ref{croce}), the map $G$ is well defined on $[0,1]\times E$, and
$G(t,\cdot)$ is a $C^k$ diffeomorphism of $E$ onto $E$ for every $t\in
[0,1]$. By (\ref{croce}),
\[
|\tilde{f}(\xi)| \leq \|T\|\, |\tilde{g}(\xi)| \leq \|T\| e^{\epsilon}
|\xi|,
\]
so the fact that $e^{\epsilon} \|T\|<1$ implies that
$\tilde{f}$ is a topological contraction of $E$, concluding the
proof.
\qed


\bibliographystyle{amsalpha}

\providecommand{\bysame}{\leavevmode\hbox to3em{\hrulefill}\thinspace}
\providecommand{\MR}{\relax\ifhmode\unskip\space\fi MR }
\providecommand{\MRhref}[2]{%
  \href{http://www.ams.org/mathscinet-getitem?mr=#1}{#2}
}
\providecommand{\href}[2]{#2}

\end{document}